\documentclass[10pt,twoside]{article}
\usepackage{latexsym}
 \usepackage{mathrsfs}
\def\thebibliography#1{\center{\bf\normalsize References}\list
 {[\arabic{enumi}]}{\settowidth\labelwidth{[#1]}\leftmargin\labelwidth
 \advance\leftmargin\labelsep
 \usecounter{enumi}}
 \def\newblock{\hskip .11em plus .33em minus .07em}
 \sloppy\clubpenalty4000\widowpenalty4000
 \sfcode`\.=1000\relax}
\def\tt{\hspace{-0.19cm}{\bf .}\hspace{0.19cm}}

\usepackage{graphicx}
\def\Proof{\noindent\hspace{2em}{\bf Proof.\hspace{1em}}}

\def\nnb{\nonumber}
\def\endpf{\hfill$\Box$\vspace{0.4cm}}

\def\ds{\displaystyle}
\newcommand{\eqref}[1]{$(\ref{#1})$}
\newcommand{\refeq}[1]{$(\ref{#1})$}
\newcommand{\thb}[1]{{\rm (#1)}}
\def\ol{\overline}
\def\GT{\Theta}

\def\cP{{\cal P}}

\def\cR{{\cal R}}
\def\cS{{\cal S}}

\def\cU{{\cal U}}

\def\cZ{{\cal Z}}
\def\mcS{{\mathscr S}}

\def\bA{\ol{A}}

\def\ba{\bar a}
\def\bb{\bar b}
\def\bc{\bar c}

\def\be{\bar e}

\def\bl{\bar l}

\def\bp{\bar p}

\def\br{\bar r}

\def\bt{\bar t}
\def\bu{\bar u}

\def\hb{\hat b}

\def\mm{\mbox{ }}

\def\bb{\hspace{18pt}}

\def\qq{\qquad}
\def\q{\quad}

\def\eqif{\mm {\rm if}\mm }
\def\eqon{ \mm {\rm on } \mm }
\def\eqin{ \mm {\rm in } \mm }

\def\all{  \mm \forall \mm }

\def\tr{\mm{\rm tr}\mm}

\def\pa{\partial}
\def\coh{\mm{\rm \overline{\, co}}\mm}

\def\defeq{\stackrel{\triangle}{=}}
\def\dbR{{\mathop{\rm l\negthinspace R}}}

\def\3n{\negthinspace \negthinspace \negthinspace }
\def\2n{\negthinspace \negthinspace }
\def\1n{\negthinspace }

\def\dbR{{\mathop{\rm l\negthinspace R}}}
\def\={\buildrel \triangle \over =}

\def\ns{\noalign{\ss}}

\def\a{\alpha}
\def\b{\beta}

\def\d{\delta}
\def\e{\varepsilon}

\def\l{\lambda}
\def\m{\mu}
\def\n{\nu}

\def\f{\varphi}
\def\th{\theta}

\def\Th{\Theta}
\def\L{\Lambda}

\def\F{\Phi}
\def\O{\Omega}

\def\cP{{\cal P}}

\def\cR{{\cal R}}
\def\cS{{\cal S}}

\def\cU{{\cal U}}

\def\cZ{{\cal Z}}

\def\ss{\smallskip}

\def\q{\quad}
\def\qq{\qquad}
\def\hb{\hbox}

\def\lan{\mathop{\langle}}
\def\ran{\mathop{\rangle}}

\def\pa{\partial}

\def\cd{\cdot}
\def\cds{\cdots}

\def\ae{\hbox{\rm a.e.{ }}}

\def\span{\hbox{\rm span$\,$}}
\def\tr{\hbox{\rm tr$\,$}}

\def\coh{\mathop{\overline{\rm co}}}

\def\div{\mathop{\nabla\cd}}

\def\({\Big (}
\def\){\Big )}
\def\[{\Big[}
\def\]{\Big]}
\def\bde{\begin{definition}}
\def\ede{\end{definition}}
\def\be{\begin{equation}}
\def\bel{\begin{equation}\label}
\def\ee{\end{equation}}
\def\bt{\begin{theorem}}
\def\et{\end{theorem}}
\def\bc{\begin{corollary}}
\def\ec{\end{corollary}}
\def\bl{\begin{lemma}}
\def\el{\end{lemma}}
\def\bp{\begin{proposition}}
\def\ep{\end{proposition}}
\def\bas{\begin{assumption}}
\def\eas{\end{assumption}}
\def\br{\begin{remark}}
\def\er{\end{remark}}
\def\ba{\begin{array}}
\def\ea{\end{array}}
\def\ed{\end{document}}

\def\square#1{\vbox{\hrule\hbox{\vrule height#1%
     \kern#1\vrule}\hrule}}
\def\rectangle#1#2{\vbox{\hrule\hbox{\vrule height#1%
     \kern#2\vrule}\hrule}}

\begin{document}
\title{
{\bf Optimality Conditions for Semilinear  Parabolic Equations with
Controls in Leading Term
\thanks{The author was supported in part by NSFC (No. 10671040 and 10831007)
and FANEDD (No. 200522).}}
 }
\author{ Hongwei Lou \thanks{School of Mathematical Sciences, and LMNS, Fudan
University, Shanghai 200433, China (Email:
\texttt{hwlou@fudan.edu.cn})}
 }

\date{}
\maketitle
\begin{quote}
\footnotesize {\bf Abstract.} An optimal control problem for
semilinear parabolic partial differential equations is considered.
The control variable appears in the leading term of the equation.
Necessary conditions for optimal controls are established by the
method of homogenizing spike variation. Results for problems with
state constraints are also stated.


\end{quote}
\normalsize

\newtheorem{Definition}{\bb Definition}[section]
\newtheorem{Theorem}[Definition]{\bb Theorem}
\newtheorem{Lemma}[Definition]{\bb Lemma}
\newtheorem{Corollary}[Definition]{\bb Corollary}
\newtheorem{Proposition}[Definition]{\bb Proposition}
\newtheorem{Remark}{\bb Remark}[section]
\newtheorem{Example}{\bb Example}
\newfont{\Bbb}{msbm10 scaled\magstephalf}
\newfont{\frak}{eufm10 scaled\magstephalf}
\newfont{\sfr}{msbm7 scaled\magstephalf}
\def\theequation{1.\arabic{equation}}
\setcounter{equation}{0} \setcounter{section}{1}
\setcounter{Definition}{0} \setcounter{Remark}{0}\textbf{1.
Introduction.}  We will give necessary conditions of optimal
controls for parabolic partial differential equation  (PDE, for
short) with leading term containing controls. This is an analogue of
the result we got for elliptic PDE with controls in the leading term
(\cite{Lou-Yong 2009}). Let us consider the following controlled  parabolic PDE
of divergence form:
\begin{equation}\label{E101}\left\{\begin{array}{ll}
\ds\pa _t z(t,x)-\div\big(A(t,x,u(t,x))\nabla
z(t,x)\big)=f(t,x,z(t,x),u(t,x)),\qq\eqin   \O_T,\\
\ds z(t,x)=0,  \qq\eqon  [0,T]\times \pa\O,\\
\ds z(0,x)=z_0(x),\qq\eqin  \O,
\end{array}\right.
\end{equation}
where $\O_T=(0,T)\times \O$, $T>0$ and $\O\subseteq\dbR^n$ is a
bounded domain with a smooth boundary $\pa\O$, $A:\O_T\times
U\to\dbR^{n\times n}$ is a map taking values in the set of all
positive definite matrices, $f:\O_T\times\dbR\times U\to\dbR$, with
$U$ being a separable metric space and $z_0(\cd)\in L^\infty(\O)$.
Function $u(\cd)$, called a {\it control}, is taken from the set
$$
\cU\equiv\big\{w:\O_T\to U\bigm|w(\cd)\hb{ is measurable }\big\}.
$$
Under some mild conditions, corresponding to a $u(\cd)\in\cU$, \eqref{E101}
admits a unique weak solution $z(\cd)\equiv z(\cd\,;u(\cd))$ which
is called the {\it state}.
We measure the performance of the control by the following cost
functional
\begin{equation}\label{E102}
J(u(\cd))=\int_{\O_T} f^0(t,x,z(t,x),u(t,x))\,
dtdx
\end{equation}
for some given map $f^0:\O_T\times\dbR\times U\to\dbR$. Our optimal
control problem is stated as follows.

{\bf Problem (C)}. Find a $\bar u(\cd)\in\cU$ such that
\begin{equation}\label{E103}
J(\bar u(\cd))=\inf_{u(\cd)\in\cU}J(u(\cd)).
\end{equation}
Any $\bar u(\cd)\in\cU$ satisfying \eqref{E103} is called an {\it
optimal control}, and the corresponding $\bar z(\cd)\equiv
z(\cd\,;\bar u(\cd))$ is called an {\it optimal state}. The pair
$(\bar z(\cd),\bar u(\cd))$ is called an {\it optimal pair}. When
$A(t,x,u)\equiv A(t,x)$, Problem (C) has been studied by many
authors, see \cite{Li-Yong 1995} and the references cited therein.
Works concerning the elliptic cases with leading term containing
controls can be founded in \cite{CCM},  \cite{Casas 1992},
\cite{C-T-J}, \cite{G-L}, \cite{Lou-Yong 2009}, \cite{M-T},
\cite{Raitums-Schmdt2005}, and \cite{Ser},  etc. However, it seems
that there are only few works devoted to parabolic cases (see
\cite{CalCas}, \cite{Tag}, etc.).

In this paper, we make the following assumptions.

(S1) Let $T>0$ and $\O$ be a bounded domain in $\dbR^n$ with a
smooth boundary $\pa\O$.

(S2) Let $U$ be a separable metric space.

(S3) Functions $A(t,x,v)=\big(a_{ij}(t,x,v)\big)$ take values in the
set $\cS^n_+$ of $ n\times n $ (symmetric) positive definite
matrices, which are measurable in $(t,x)\in \O_T$ and continuous in
$v\in U$. Moreover, there exist $\L\geq\l>0$ such that for almost all
$(t,x)\in \O_T$,
\begin{equation}\label{E104}\begin{array}{ll}
\ds\l|\xi|^2\le\lan A(t,x,v)\xi,\xi\ran\leq\L |\xi|^2,\qq\forall
\xi\in\dbR^n,  \,v\in U, \end{array}
\end{equation}
where $\lan\cd\,,\cd\ran$ stands for the inner product in $\dbR^n$.

(S4) Functions $f(t,x,z,v)$ is measurable in $(t,x)$,   continuous
in $(z,v)\in\dbR\times U$, and continuously differentiable in $z$.
Moreover,   there exits a constant $M>0$ such that
\begin{equation}\label{E105}
zf(t,x,z,v)\le M(z^2+1),\qq\forall (t,x,z,v)\in
\O_T\times\dbR\times U
\end{equation}
and for any $R>0$, there exists an $M_R>0$ such that
\begin{equation}\label{E106}
|f(t,x,z,v)|+|f_z(t,x,z,v)|\leq M_R, \qq\ae (t,x,v)\in
       \O_T\times U, \, |z|\le R.
\end{equation}

(S5) Function $f^0(t,x,z,v)$ is measurable in $(t,x)$,  continuous
in $(z,v)\in\dbR\times U$, and continuously differentiable in $z$.
Moreover,  for any $R>0$, there exists a $K_R>0$ such that
\begin{equation}\label{E107}
|f^0(t,x,z,v)|+|f^0_z(t,x,z,u)|\le K_R,\qq\ae(t,x,v)\in
       \O_T\times U, \, |z|\le R.
\end{equation}

Our main result is the following.

\begin{Theorem}\label{T101}\tt Let {\rm(S1)--(S5)} hold and $z_0\in
L^\infty(\O)$. Let $(\bar z(\cd),\bar u(\cd))$ be an optimal pair of
Problem {\rm(C)}. Let $\bar\psi(\cd)$ be the weak solution of the
following adjoint equation
\begin{equation}\label{adjoint}
\left\{\begin{array}{l} \ds \pa_t \bar \psi (t,x)  \ds
+\div\big(A(t,x,\bar u(t,x))\nabla\bar\psi(t,x)\big)=  \ds
f^0_z(t,x,\bar z(t,x),\bar
u(t,x))\\
  \ds \hspace{4em} -f_z(t,x,\bar z(t,x),\bar u(t,x))\,\bar\psi(t,x)
,\qq\eqin  \O_T,\\
\ds\bar\psi(t,x) =0,\qq\eqon  [0,T]\times \pa\O, \\
\ds \bar\psi(T,x)=0, \qq\eqin   \O.
\end{array}\right.
\end{equation}
Then
\begin{eqnarray}\label{E109}
\nnb & \ds H\big(t,x,\bar z(t,x),\bar\psi(t,x),\nabla\bar
z(t,x),\nabla\bar\psi(t,x),\bar u(t,x)\big)
-H\big(t,x,\bar z(t,x),\bar\psi(t,x),\nabla\bar z(t,x),\nabla\bar\psi(t,x),v\big)\\
\nnb\ds\ge & \ds {1\over2}\big|A(t,x,v)^{-{1\over 2}}(A(t,x,\bar
u(t,x))-A(t,x,v))\nabla\bar z(t,x)\big|\,\big|A(t,x,v)^{-{1\over
2}}(A(t,x,\bar u(t,x))-A(t,x,v))\nabla\bar\psi(t,x)\big|\\
\nnb & \ds +{1\over2}\lan A(t,x,v)^{-{1\over2}}(A(t,x,\bar
u(t,x))-A(t,x,v))\nabla\bar z(t,x),A(t,x,v)^{-{1\over2}}(A(t,x,\bar
u(t,x))-A(t,x,v))\nabla\bar\psi(t,x)\ran,\\
 & \ds\qq\qq\qq \forall v\in U,\q\ae
(t,x)\in\O_T,
\end{eqnarray}
where
\begin{eqnarray}\label{E110}
\nnb && \ds H(t,x,z,\psi,\xi,\eta,v)
=\lan\psi,f(t,x,z,v)\ran-f^0(t,x,z,v)
-\lan A(t,x,v)\xi,\eta\ran,\\
&& \ds   \qq\qq\qq (t,x,z,\psi,\xi,\eta, v)\in [0,T]\times
\O\times\dbR\times \dbR\times\dbR^n\times\dbR^n\times U.
\end{eqnarray}
\end{Theorem}
Since the right hand side of (\ref{E109}) is always nonnegative,
\eqref{E109} implies
\begin{eqnarray}\label{E111}
\nnb & &\ds H\big(t,x,\bar z(t,x),\bar\psi(t,x),\nabla\bar
z(t,x),\nabla\bar\psi(t,x),\bar u(t,x)\big) \\
\nnb &= &\ds \max_{v\in U} H\big(t,x,\bar
z(t,x),\bar\psi(t,x),\nabla\bar
z(t,x),\nabla\bar\psi(t,x),v\big), \\
&& \qq\qq\qq\qq\forall v\in U,\q\ae (t,x)\in \O_T.
\end{eqnarray}
When $A(t,x,v)\equiv A(t,x)$, the right hand side of (\ref{E109})
 is zero, thus, the result automatically recovers those for the classical semilinear
case without state constraints (\cite{Li-Yong 1995}).

Since $U$ is not necessarily convex, it is well-known that people usually
use spike variations to derive necessary conditions for optimal
controls. Such a spike variation technique does not directly work
for problems with leading term containing the control. To overcome
the difficulty, we adopt the idea of homogenization for  PDEs to
carefully select some special type spike variations of controls so
that we can have desired ``differentiability'' of the state with
respect to the control. We  can see in \cite{Lou-Yong 2009} that
such a method is useful
 for the cases of elliptic PDEs. The
main idea to treat parabolic case is same to that for elliptic case.
However, there are some new difficulties in studying properties of
variational equations.

Comparing Theorem \ref{T101} and the corresponding result for
elliptic case in \cite{Lou-Yong 2009}, we can see that they are
similar when $n\geq 2$ and  slightly different when $n=1$. More
precisely, Theorem \ref{T101} of this paper is very similar to
Theorem 1.1 in \cite{Lou-Yong 2009} for high dimensional cases. In
particular, for parabolic case with $n=1$, instead of
\begin{eqnarray}\label{E112}
\nnb & & \ds H\big(t,x,\bar z(t,x),\bar\psi(t,x), \bar z_x (t,x),
\bar\psi_x(t,x),\bar u(t,x)\big)
-H\big(t,x,\bar z(t,x),\bar\psi(t,x), \bar z_x(t,x), \bar\psi_x(t,x),v\big)\\
 \ds & \ge & \ds  { (A(t,x,\bar u(t,x))-A(t,x,v))^2 \over A(t,x,v)}\bar z_x(t,x)
\bar\psi_x(t,x),  \qq \all v\in U,\q\ae (t,x)\in\O_T,
\end{eqnarray}
we have \refeq{E109}, i.e.,
\begin{eqnarray}\label{E113}
\nnb & & \ds H\big(t,x,\bar z(t,x),\bar\psi(t,x), \bar z_x (t,x),
\bar\psi_x(t,x),\bar u(t,x)\big)
-H\big(t,x,\bar z(t,x),\bar\psi(t,x), \bar z_x(t,x), \bar\psi_x(t,x),v\big)\\
 \ds & \ge & \ds   { (A(t,x,\bar u(t,x))-A(t,x,v))^2 \over
A(t,x,v)}\Big[\bar z_x(t,x) \bar\psi_x(t,x)\Big]^+,  \qq \all v\in
U,\q\ae (t,x)\in\O_T.
\end{eqnarray}
One can see that \refeq{E112} is similar to the corresponding result
for elliptic case with $n=1$, while \refeq{E109} (i.e.,
\refeq{E113}) is similar to the corresponding result for elliptic
case with $n\geq 2$. We mention that for elliptic cases with $n\geq
2$, the corresponding right hand of \refeq{E109} follows from a fact
given in Lemma \ref{T205}. While for parabolic case with $n=1$, the
right hand of \refeq{E109} (i.e.,  \refeq{E113}) follows in a
different way. In fact, it follows from \refeq{E112} and
\begin{eqnarray}\label{E114}
\nnb & & \ds H\big(t,x,\bar z(t,x),\bar\psi(t,x), \bar z_x (t,x),
\bar\psi_x(t,x),\bar u(t,x)\big)
-H\big(t,x,\bar z(t,x),\bar\psi(t,x), \bar z_x(t,x), \bar\psi_x(t,x),v\big)\\
 \ds & \ge & \ds  0,  \qq\qq \all v\in
U,\q\ae (t,x)\in\O_T.
\end{eqnarray}
From the proof of Theorem \ref{T101}, one can see that \refeq{E112}
can be yielded from using spike variation along space-direction and
\refeq{E114} can be yielded from using spike variation along
time-direction (see \refeq{E317}).

Another difference between parabolic cases and elliptic cases appear
in that there are three possible types of homogenized equations for
parabolic cases when taking a different scale for the time and the
space variables, while there is only one type of homogenized
equations for elliptic cases. Difficulty occurs in analyzing the
second type of homogenized equations (see the proof of Lemma
\ref{T202} for details). Despite the different types of homogenized
equations, the variational equations are same and we finally get
same  optimality conditions for the three cases.  Nevertheless, we
think results of this paper will be useful to analyze  the
second-order variational equations,
which is a problem more difficult than that for  first-order
variational equations.

The rest of the paper is organized as follows. In Section 2, we
present some preliminary results. Section 3 is devoted to a proof of
our main result. Problem with state constraints will be
discussed in Section 4. 

\vspace{6mm}

\def\theequation{2.\arabic{equation}}
\setcounter{equation}{0} \setcounter{section}{2}
\setcounter{Definition}{0} \setcounter{Remark}{0}\textbf{2.
Preliminaries.} In this section, we will give some preliminary
results needed in proving Theorem \ref{T101}.  For $Y=[0,\a_1]\times
[0,\a_2]\times \ldots [0,\a_n]$, a function $g(x)$ on $\dbR^n$ is
called $Y$-periodic if it admits period $\a_j$ in the direction
$x_j$ ($j=1,2,\ldots,n$).

\begin{Lemma}\tt\label{T201} Let $r>0$, $\d\in (0,1)$ and
 {\rm(S1)} hold. Let $h(\cd)\in L^2(\O_T)$, and
$A^m(\cd)=\big(a^m_{ij}(\cd)\big)\in L^\infty(\O_T;\cS^n_+)$  such
that for some $\L\ge\l>0$,
\begin{equation}\label{E201}
\l|\xi|^2\le\lan A^m(t,x)\xi,\xi\ran
\le\L|\xi|^2,\qq\forall \xi\in \dbR^n,(t,x)\in \O_T,
m=1,2,3,4.
\end{equation}
Define
\begin{equation}\label{FE202} G(t,x,s,y)\equiv \big(g_{ij}(t,x,s,y)\big)\equiv
\big(g_{ij}(t,x,s,y_1)\big)=\left\{\begin{array}{ll}\ds
A^1(t,x),& \eqif  (\{s\},\{y_1\})\in [\d,1)\times [\d,1),\\
\ds A^2(t,x),& \eqif   (\{s\},\{y_1\})\in [0,\d)\times [\d,1),\\
\ds A^3(t,x),& \eqif   (\{s\},\{y_1\})\in [\d,1)\times [0,\d),\\
\ds A^4(t,x),& \eqif   (\{s\},\{y_1\})\in [0,\d)\times
[0,\d),\end{array}\right.
\end{equation}
where $\{a\}$ denote the decimal part of a real number $a$. For
$\e>0$, let $z^\e(\cd)\in L^2(0,T;H^1_0(\O))$ be the weak solution
of
\begin{equation}\label{E202}\left\{\begin{array}{ll}
\ds\pa_t z^\e(t,x)-\div\big[G \big(t,x,{t\over \e^r},{x\over
\e}\big)\nabla
z^\e(t,x)\big]=h(t,x),\qq\eqin \O_T,\\
\ds z^\e(t,x)=0,  \qq\eqon
 [0,T]\times \pa\O,\\
 \ds z^\e (0,x)=z_0(x),\qq\eqin  \O
\end{array}\right.
\end{equation}
with $z_0(\cd)\in L^\infty(\O)$.  Then
\begin{equation}\label{E203}
z^\e(\cd)\to z(\cd),\qq\hb{weakly in}\q
L^2(0,T;H^1_0(\O))
 \end{equation}
with $z(\cd)$ being the weak solution of
\begin{equation}\label{E204}\left\{\begin{array}{ll}
\ds\pa_t z (t,x)-\div\big(Q (t,x)\nabla
z (t,x)\big)=h(t,x),\qq\eqin  \O_T,\\
\ds z (t,x)=0,  \qq\eqon  [0,T]\times \pa\O,\\
\ds z  (0,x)=z_0(x),\qq\eqin  \O,
\end{array}\right.
\end{equation}
and $Q(\cd)=\big(q_{ij}(\cd)\big)\in L^\infty(\O_T;\cS^n_+)$ being
given by
\begin{equation}\label{E205}q_{ij}(t,x)=\int^1_0dy_1\int^1_0 \big(g_{ij}(t,x,s,y_1)+
g_{i1}(t,x,s,y_1)\pa_{y_1}\f^j(t,x,s,y_1)\big)ds \,,\qq 1\le i,j\le
n,
\end{equation}
where $\f^k(t,x,\cd)\in
L_\#^2(0,1;W^{1,2}_\#(0,1)/\dbR)$, $\#$ means the function is
$[0,1]$ periodic.

For $r<2$,  $\f^k(t,x,\cd)$ is the unique solution of
\begin{equation}\label{E206} \pa_{y_1}
\big(g_{1k}(t,x,s,y_1)+g_{11}(t,x,s,y_1)\pa_{y_1} \f^k
(t,x,s,y_1)\big)=0.
\end{equation}

For $r=2$,  $\f^k(t,x,\cd)$ is the solution of
\begin{equation}\label{E207}
\pa_s \f^k(t,x,s,y_1)-\pa_{y_1}
\big(g_{1k}(t,x,s,y_1)+g_{11}(t,x,s,y_1)\pa_{y_1} \f^k
(t,x,s,y_1)\big)=0.
\end{equation}

While for $r>2$, $\f^k(t,x,s,y)=\f^k(t,x,y_1)$ is the solution of
\begin{equation}\label{E208}
\pa_{y_1} \Big(\int^1_0g_{1k}(t,x,s,y_1)\, ds +\int^1_0
g_{ij}(t,x,s,y_1)\, ds \, \pa_{y_1} \f^k
(t,x,y_1)\Big)=0.
\end{equation}
\end{Lemma}
\Proof The above proposition is a corollary of Theorem 2.1 in
\cite{B-L-P}, Chapter 2 (see also Remark 1.1 and ``Comments and
Problems" there). The result can also be got by the technique of two
scale convergence (\cite{Alla}, \cite{Hol}).   According to
\cite{B-L-P},
\begin{equation}\label{E209}
q_{ij}(t,x)=\int_{[0,1]^n}dy\int^1_0 \big(g_{ij}(t,x,s,y)+
\sum^n_{k=1}g_{ik}(t,x,s,y)\pa_{y_k}\f^j(t,x,s,y)\big)ds \,,\qq 1\le
i,j\le n.
\end{equation}

For $r<2$, $\f^k(t,x,\cd)\in L_\#^2(0,1;W^{1,2}_\#((0,1)^n)/\dbR)$
is the unique solution of
\begin{equation}\label{E210} \sum^n_{i,j=1}\pa_{y_i}
\big(g_{ij}(t,x,s,y)\d_{jk}+g_{ij}(t,x,s,y)\pa_{y_j} \f^k
(t,x,s,y)\big)=0,
\end{equation}
where $\d_{ij}$ equals to $1$ if $i=j$ and $0$ if $i\ne j$.

For $r=2$,  $\f^k(t,x,\cd)$ is the solution of
\begin{equation}\label{E211}
\pa_s \f^k(t,x,s,y)-\sum^n_{i,j=1}\pa_{y_i}
\big(g_{ij}(t,x,s,y)\d_{jk}+g_{ij}(t,x,s,y)\pa_{y_j} \f^k
(t,x,s,y)\big)=0.
\end{equation}

While for $r>2$, $\f^k(t,x,s,y)=\f(t,x,y)$ is the solution of
\begin{equation}\label{E212}
 \sum^n_{i,j=1}\pa_{y_i}
\Big(\int^1_0g_{ij}(t,x,s,y)\d_{jk}\, ds +\int^1_0 g_{ij}(t,x,s,y)\,
ds \, \pa_{y_j} \f^k (t,x,y)\Big)=0.
\end{equation}

Since $G(t,x,s,y)$ is independent of $y_2,y_3,\ldots, y_n$, we must
have $\f^k(t,x,s,y)=\f^k(t,x,s,y_1)$ and consequently,
(\ref{E209})---(\ref{E212}) becomes (\ref{E205})---(\ref{E208}).
\endpf

The following lemma  concerns the ``derivative" of $q_{ij}$ in
$\d=0$.

\begin{Lemma}\tt\label{T202} Let $r>0$, $\L>\l>0$. Assume $\L\ge a_m\ge\l $,
$|b_m|\le \L$, $|c_m|\le \L$ \thb{$m=1,2,3,4$}. Let $\d\in (0,1)$
and
\begin{equation}\label{} (a^\d(s,y),b^\d(s,y),c^\d(s,y))=\left\{\begin{array}{ll}\ds
(a_1,b_1,c_1),& \eqif  (\{s\},\{y\})\in [\d,1)\times [\d,1),\\
\ds (a_2,b_2,c_2),& \eqif   (\{s\},\{y\})\in [0,\d)\times [\d,1),\\
\ds (a_3,b_3,c_3),& \eqif   (\{s\},\{y\})\in [\d,1)\times [0,\d),\\
\ds (a_4,b_4,c_4),& \eqif   (\{s\},\{y\})\in [0,\d)\times
[0,\d).\end{array}\right.
\end{equation}
\begin{center}
\includegraphics[height=5.0cm]{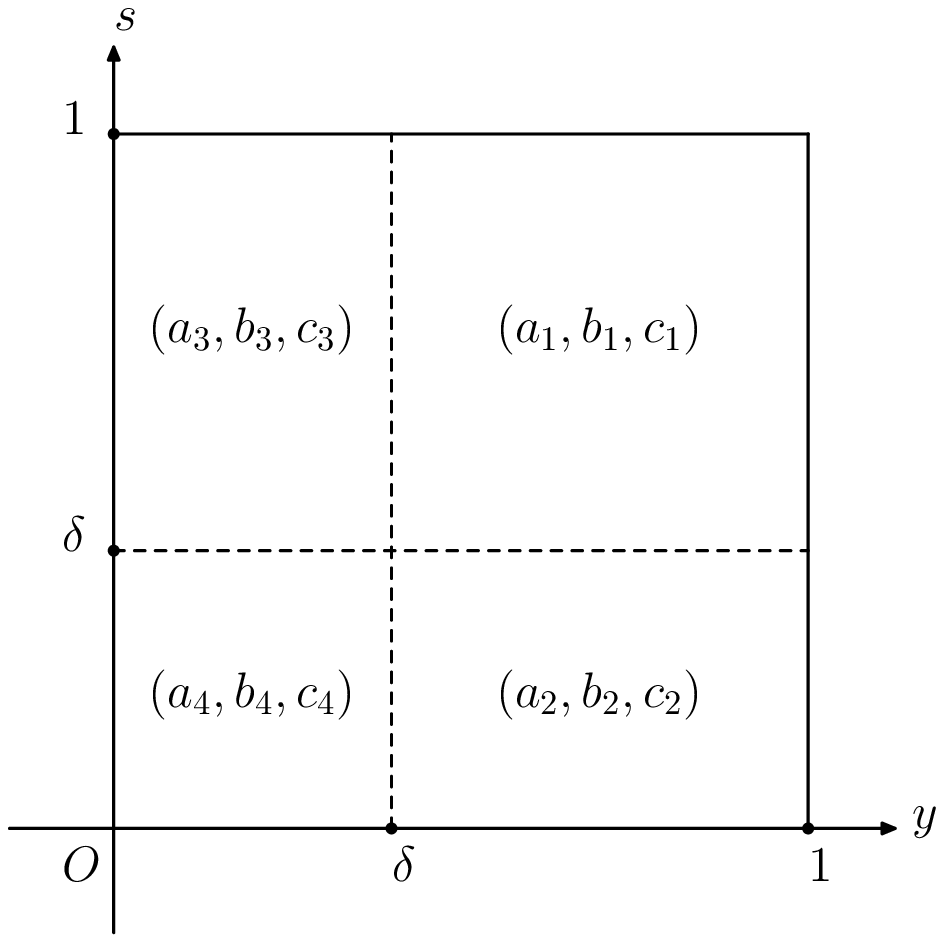}
\end{center}
Let  $\phi_1^\d(\cd)\in L_\#^2(0,1;W^{1,2}_\#(0,1)/\dbR)$ be the
solution of
\begin{equation}\label{E213} \pa_y \big(b^\d(s,y)+a^\d(s,y)\pa_y \phi_1^\d
(s,y)\big)=0,
\end{equation}
$\phi_2^\d(\cd)$ be the solution of
\begin{equation}\label{E214}
 \pa_s \phi_2^\d(s,y)-\pa_y
\big(b^\d(s,y)+a^\d(s,y)\pa_y\phi_2^\d (s,y)\big)=0
\end{equation}
and $\phi_3^\d(s,y)\equiv \phi_3^\d(y)$ be the solution of
\begin{equation}\label{E215}
\pa_y \Big(\int^1_0b^\d(s,y)\, ds +\int^1_0 a^\d
(s,y)\, ds \, \pa_y \phi_3^\d (y)\Big)=0.
\end{equation}
Then there exists a constant $C=C(\L,\l)$, independent of $\d\in
(0,1)$, such that
\begin{equation}\label{E216}
 \Big|{1\over \d}\int^1_0dy\int^1_0 c^\d(s,y)\pa_y
\phi_k^\d(s,y)\, ds+{(b_3-b_1)(c_3-c_1)\over a_3}\Big|\leq C\sqrt
\d,\qq k=1,2,3.
\end{equation}
\end{Lemma}
\Proof.  \textbf{I.}  It follows from (\ref{E213}) that
$$
b^\d(s,y)+a^\d(s,y)\pa_y \phi_1^\d (s,y)=p^\d(s).
$$
Since $\phi_1^\d(\cd)$ is $[0,1]^2$ periodic, we get that
$$
\int^1_0 {p^\d(s)-b^\d(s,y)\over a^\d(s,y)}\, dy=\int^1_0 \pa_y
\phi_1^\d (s,y)\, dy=0,\qq s\in [0,1].
$$
Solve $p^\d(\cd)$ from the above we  get
\begin{equation}\label{E217}
 \pa_y \phi_1(s,y)= \left\{\begin{array}{ll}\ds {b_3-b_1\over
(1-\d)a_3+\d a_1}\d, & (s,y)\in (\d,1)\times (\d,1),\\
\ds {b_4-b_2\over
(1-\d)a_4+\d a_2}\d, & (s,y)\in (0,\d)\times (\d,1),\\
\ds -{b_3-b_1\over
(1-\d)a_3+\d a_1}(1-\d), & (s,y)\in (\d,1)\times (0,\d),\\
\ds -{b_4-b_2\over (1-\d)a_4+\d a_2}(1-\d), & (s,y)\in (0,\d)\times
(0,\d).\end{array}\right.
\end{equation}
Thus
\begin{equation}\label{E218} |\pa_y \phi_1^\d(s,y)|\le \left\{\begin{array}{ll}\ds {2\L\over
\l}\d,   & (s,y)\in [0,1]\times (\d,1),\\
 & \\
\ds {2\L\over \l}(1-\d), & (s,y)\in [0,1]\times
(0,\d).\end{array}\right.
\end{equation}
A direct calculation shows that
\begin{eqnarray}\label{E219}
\nnb  && \ds   \Big|{1\over \d}\int^1_0dy\int^1_0 c^\d(s,y)\pa_y
\phi_1^\d(s,y)\, ds+{(b_3-b_1)(c_3-c_1)\over
a_3}\Big|\\
\nnb &=& \ds \Big|\d (2-\d){(b_3-b_1)(c_3-c_1)\over (1-\d) a_3+\d
a_1}-\d (1-\d){(b_4-b_2)(c_4-c_2)\over (1-\d) a_4+\d
a_2}\Big|\\
  &\le &\ds \d(2-\d){4\L^2\over \l}+\d(1-\d) {4\L^2\over \l}\le
{12\L^2\over \l}\d.
\end{eqnarray}

\textbf{II.} It follows from (\ref{E214}) and the periodicity of
$\phi^\d_2(\cd)$ that
\begin{eqnarray}\label{E220}
\nnb 0&=& \ds \int^1_0\int^1_0 [
b^\d(s,y)\pa_y\phi_2^\d(s,y) +a^\d(s,y)|\pa_y\phi_2^\d(s,y)|^2]\,
ds\,  dy\\
&=&\ds \int^1_0\int^1_0 [ (b^\d(s,y)-b_1)\pa_y\phi_2^\d(s,y)
+a^\d(s,y)|\pa_y\phi_2^\d(s,y)|^2]\, ds\,  dy.
\end{eqnarray}
Thus,
\begin{eqnarray}\label{}  \nnb && \ds \int^1_0\int^1_0  |\pa_y\phi_2^\d(s,y)|^2\, \,
ds\,  dy\\
\nnb &\le & \ds {1\over \l}\int^1_0\int^1_0
a^\d(s,y)|\pa_y\phi_2^\d(s,y)|^2\, \, ds\, dy\\
\nnb&= &\ds -{1\over \l}\int^1_0\int^1_0
(b^\d(s,y)-b_1)\pa_y\phi_2^\d(s,y)
\, ds\,  dy\\
\nnb &\le &\ds {1\over \l}\Big(\int^1_0\int^1_0 |b^\d(s,y)-b_1|^2\,
ds\, dy\Big)^{1\over 2}\, \Big(\int^1_0\int^1_0
|\pa_y\phi_2^\d(s,y)|^2\, ds\,
dy\Big)^{1\over 2}\\
& \le &\ds {2\sqrt 2 \L \sqrt \d\over \l} \, \Big(\int^1_0\int^1_0
|\pa_y\phi_2^\d(s,y)|^2\, ds\, dy\Big)^{1\over 2}.
\end{eqnarray}
Therefore,
\begin{equation}\label{}
\|\pa_y\phi_2^\d\|_{L^2([0,1]^2)} \le   {2\sqrt 2 \L \sqrt
\d\over \l}.
\end{equation}
On the other hand, \eqref{E218} implies
\begin{equation}\label{}
\|\pa_y\phi_1^\d\|_{L^2([0,1]^2)} \le   {2\L \sqrt \d\over
\l}.
\end{equation}
Denote $\Phi^\d(\cd)=\phi_2^\d(\cd)-\phi_1^\d(\cd)$, we have
\begin{equation}\label{}  \|\pa_y\Phi^\d\|_{L^2([0,1]^2)} \le   {2(1+\sqrt 2)\L \sqrt
\d\over \l}.
\end{equation}
For any $\f(\cd)\in W^{1,2}_\#(0,1)$, it follows from
(\ref{E213})---(\ref{E214}) that
\begin{equation}\label{} \int^1_0  [b^\d(s,y)+a^\d(s,y)\pa_y\phi_1^\d(s,y)]\, \pa_y \f
(y)\, dy=0, \qq \forall s\in [0,1],
\end{equation}
\begin{equation}\label{} \int^1_0  \int^1_0 [b^\d(s,y)+a^\d(s,y)\pa_y\phi_2^\d(s,y)]\,
\pa_y \f (y)\, ds\, dy=0.
\end{equation}
Thus,
\begin{equation}\label{} \int^1_0  \int^1_0   a^\d(s,y)\pa_y\Phi^\d(s,y) \, \pa_y \f
(y)\, ds\, dy=0.
\end{equation}
Therefore $\ds  \int^1_0   a^\d(s,y)\pa_y\Phi^\d(s,y) \, ds $ is a
constant and consequently,
\begin{eqnarray}\label{E228}
\nnb && \ds \Big|\int^1_0   a^\d(s,y)\pa_y\Phi^\d(s,y) \, ds \Big|\\
\nnb &=&\ds
\Big|\int^1_0\int^1_0   a^\d(s,w)\pa_y\Phi^\d(s,x) \, ds \, dx\Big|\\
\nnb &=&\ds
\Big|\int^1_0\int^1_0   (a^\d(s,x)-a_1)\pa_y\Phi^\d(s,x) \, ds \, dx\Big|\\
\nnb &\le &\ds \Big(\int^1_0\int^1_0   |a^\d(s,x)-a_1|^2\, ds\,
dx\Big)^{1\over 2} \, \|\pa_y\Phi^\d\|_{L^2([0,1]^2)}\\
&\le &\ds \sqrt 2\L\sqrt \d \cd {2(1+\sqrt 2)\L\over \l}\sqrt
\d={2(2+\sqrt 2)\L^2\over \l}\d,\qq \forall y\in
[0,1].
\end{eqnarray}
Nothing that (\ref{E213}) implies
\begin{equation}\label{}   \int^1_0  [ b^\d(s,y)\pa_y\phi_2^\d(s,y) +a^\d(s,y)
\pa_y\phi_1^\d(s,y)\pa_y\phi_2^\d(s,y) ]\,
 dy=0, \qq\forall s\in [0,1],
 \end{equation}
we get from (\ref{E220}) that
$$
\int^1_0\int^1_0 a^\d(s,y) \pa_y\Phi^\d(s,y)\pa_y\phi_2^\d(s,y)  \,
 ds\, dy=0.
$$
Thus denote
$$
\th_{1,\d}={b_3-b_1\over (1-\d)a_3+\d a_1},
\q\th_{2,\d}={b_4-b_2\over (1-\d)a_4+\d a_2},
$$
we have
\begin{eqnarray*}
 & & \ds \|\pa_y\Phi^\d\|_{L^2([0,1]^2;a^\d)}^2 \equiv    \ds  \int^1_0\int^1_0
a^\d(s,y)|\pa_y\Phi^\d(s,y)|^2\, \, ds\, dy\\
&= &\ds -\int^1_0\int^1_0 a^\d(s,y)
\pa_y\Phi^\d(s,y)\pa_y\phi_1^\d(s,y)  \,
 ds\, dy\\
&= &\ds -\int^1_\d\int^1_\d a^\d(s,y)
\pa_y\Phi^\d(s,y)\th_{1,\d}\d\, ds \,dy -\int^1_\d\int^\d_0
a^\d(s,y) \pa_y\Phi^\d(s,y)\th_{2,\d}\d\,
ds \,dy\\
&&\ds +
   \int^\d_0\int^1_0 a^\d(s,y) \pa_y\Phi^\d(s,y)\th_{1,\d}(1-\d)\,
ds \,dy \\
&&\ds +\int^\d_0\int^\d_0 a^\d(s,y)
\pa_y\Phi^\d(s,y)(\th_{2,\d}-\th_{1,\d})(1-\d)\,
ds \,dy\\
&\le & \ds  \int^1_\d\int^1_0 a^\d(s,y) |\pa_y\Phi^\d(s,y)|\,\cd
{2\L\over \l}\d\, ds \,dy  +
   \int^\d_0\Big|\int^1_0 a^\d(s,y) \pa_y\Phi^\d(s,y)\, ds\Big| \, \cd {2\L\over \l}(1-\d)\,
dy \\
&&\ds +\int^\d_0\int^\d_0 a^\d(s,y) |\pa_y\Phi^\d(s,y)|\,   \cd
{4\L\over \l} (1-\d)\,
ds \,dy\\
&\le &\ds {2\L\over \l}\d\|\pa_y\Phi^\d\|_{L^2([0,1]^2;a^\d)}\cd
\sqrt \L +
 {2(2+\sqrt 2)\L^2\over \l}\d\, \cd {2\L\over \l}\,\d \\
&&\ds+ {4\L\over \l}\, \cd
 \sqrt \L \,\d \Big(\int^\d_0\int^\d_0 a^\d(s,y) |\pa_y\Phi^\d(s,y)|^2\,ds\, dy\Big)^{1\over 2}
\\
&\le &\ds {6\L\sqrt \L\over
\l}\,\d\|\pa_y\Phi^\d\|_{L^2([0,1]^2;a^\d)}+{4(2+\sqrt 2)\L^3\over
\l^2}\, \d^2.
\end{eqnarray*}
Thus
\begin{equation}\label{}
  \|\pa_y\Phi^\d\|_{L^2([0,1]^2)}\le {1\over \sqrt \l}\|\pa_y\Phi^\d \|_{L^2([0,1]^2;a^\d)}^2
\le {8\L\sqrt \L\over \l\sqrt \l}\, \d.
\end{equation}
Then, (\ref{E228}) can be improved as
\begin{equation}\label{}  \ds \Big|\int^1_0   a^\d(s,y)\pa_y\Phi^\d(s,y) \, ds \Big|
\le  \sqrt 2\L\sqrt \d \cd {8\L\sqrt \L\over \l \sqrt
\l}\d={12\L^2\sqrt \L\over \l\sqrt \l }\d\sqrt \d,\qq \forall y\in
[0,1].
\end{equation}
Moreover,
\begin{eqnarray}\label{}
\nnb &&\ds \Big(\int^1_0\Big|\int^1_0
\pa_y\Phi^\d(s,y) \, ds \Big|^2\, dy\Big)^{1\over 2}\\
 \nnb &\le &\ds {1\over \l}\Big(\int^1_0\Big|\int^1_0 \Big(a_1\chi_{[\d,1]}(y)+a_3\chi_{(0,\d)}(y)\Big)
\pa_y\Phi^\d(s,y) \, ds \Big|^2\, dy\Big)^{1\over 2}
\\
\nnb &\le &\ds
  {1\over \l}\Big(\int^1_0\Big|\int^1_0 a^\d(s,y)\pa_y\Phi^\d(s,y) \, ds \Big|^2\, dy\Big)^{1\over 2}\\
\nnb & &\ds  + {1\over \l}
 \Big(\int^1_0\Big|\int^\d_0 \Big((a_2-a_1)\chi_{[\d,1]}(y)+(a_4-a_3)\chi_{(0,\d)}(y)\Big)
\pa_y\Phi^\d(s,y) \, ds \Big|^2\, dy\Big)^{1\over 2}
\\
\nnb &\le &\ds {12\L^2\sqrt \L\over \l^2\sqrt \l }\d\sqrt \d+
{1\over \l} \Big[\int^1_0 \Big(4\L^2\d \int^\d_0 |\pa_y\Phi^\d(s,y)
|^2\, ds \Big)\,
dy\Big]^{1\over 2}\\
  &\le &\ds {12\L^2\sqrt \L\over \l^2\sqrt \l }\d\sqrt \d+
{2\L\sqrt \d\over \l} \cd {8\L\sqrt \L\over \l\sqrt \l}\, \d
={28\L^2\sqrt \L\over \l^2\sqrt \l }\d\sqrt \d.
\end{eqnarray}
Therefore,
\begin{eqnarray}\label{E233}
\nnb & & \ds \Big|{1\over \d}\int^1_0\int^1_0 c^\d(s,y)\pa_y
\phi_2^\d(s,y)\, ds\, dy+{(b_3-b_1)(c_3-c_1)\over
a_3}\Big|\\
\nnb &\le & \ds \Big|{1\over \d}\int^1_0dy\int^1_0 c^\d(s,y)\pa_y
\phi_1^\d(s,y)\, ds+{(b_3-b_1)(c_3-c_1)\over a_3}\Big|\\
\nnb &&\ds  +
 \Big|{1\over \d}\int^1_0\int^1_0 c^\d(s,y)\pa_y
\Phi^\d(s,y)\, ds\, dy\Big|\\
\nnb &\le & \ds {12\L^2\over \l}\,\d+\Big|{1\over
\d}\int^1_0\int^1_0
\big(c^\d(s,y)-c_1\chi_{[\d,1]}(y)-c_3\chi_{(0,\d)}(y)\big)\pa_y
\Phi^\d(s,y)\, ds\, dy\Big| \\
\nnb &&+\ds\Big|{1\over \d}\int^1_0\int^1_0
\big(c_1\chi_{[\d,1]}(y)+c_3\chi_{(0,\d)}(y)\big)\pa_y
\Phi^\d(s,y)\, ds\, dy\Big| \\
\nnb &=&\ds {12\L^2\over \l}\,\d+\Big|{1\over \d}\int^1_0\int^\d_0
\big((c_2-c_1)\chi_{[\d,1]}(y)+(c_4-c_3)\chi_{(0,\d)}(y)\big)\pa_y
\Phi^\d(s,y)\, ds\, dy\Big| \\
\nnb &&+\ds\Big|{1\over \d}\int^1_0\Big(\int^1_0 \pa_y
\Phi^\d(s,y)\, ds\Big)\, \big(c_1\chi_{[\d,1]}(y)+c_3\chi_{(0,\d)}(y)\big) \, dy\Big| \\
\nnb &\le &\ds {12\L^2\over \l}\,\d+{2\L\sqrt \d\over \d}
\|\pa_y\Phi^\d\|_{L^2([0,1]^2)}+{\L\over \d} \Big(\int^1_0\Big|
\int^1_0 \pa_y
\Phi^\d(s,y)\, ds\Big|^2\, dy\Big)^{1\over 2}\\
\nnb &\le &\ds {12\L^2\over \l}\,\d+{2\L\sqrt \d\over \d}\, \cd
{8\L\sqrt \L\over \l\sqrt\l}\,\d+  {\L\over \d}\cd {28\L^2\sqrt
\L\over \l^2\sqrt \l }\d \sqrt
\d\\
 &\le &\ds {56\L^3\sqrt \L\over \l^2\sqrt \l } \sqrt \d.
\end{eqnarray}

\textbf{III.} By (\ref{E215}), $$ \int^1_0b^\d(s,y)\, ds +
\int^1_0a^\d(s,y)\, ds\,\pa_y \phi_3^\d(y)$$ is a constant. Then
similar to (\ref{E217}), we can get from the periodicity of
$\phi_3^\d(\cd)$ that
\begin{equation}\label{}  \pa_y \phi^\d_3(y)= \left\{\begin{array}{ll}\ds
{(b_3-b_1)(1-\d)+(b_4-b_2)\d\over \d (1-\d)a_1+\d^2
a_2+(1-\d)^2a_3+\d(1-\d)a_4}\d, & y\in (\d,1),\\ \, &  \\
\ds -{(b_3-b_1)(1-\d)+(b_4-b_2)\d\over \d (1-\d)a_1+\d^2
a_2+(1-\d)^2a_3+\d(1-\d)a_4}(1-\d), & y \in
(0,\d).\end{array}\right.
\end{equation}
Thus
$$
\begin{array}{ll} & \ds  {1\over \d}\int^1_0\int^1_0 c^\d(s,y)\pa_y
\phi_3^\d(s,y)\, ds\, dy\\
 =&\ds -{(1-\d)\Big[ (1-\d) (b_3-b_1)+\d (b_4-b_2)\Big] \, \cd
\Big[ (1-\d) (c_3-c_1)+\d (c_4-c_2)\Big]\over \d (1-\d)a_1+\d^2
a_2+(1-\d)^2a_3+\d(1-\d)a_4} .\end{array}
$$
One can verify that
\begin{eqnarray}\label{E235}
\nnb & & \ds  \Big|{1\over \d}\int^1_0\int^1_0 c^\d(s,y)\pa_y
\phi_3^\d(s,y)\, ds\, dy+{  (b_3-b_1) (c_3-c_1) \over
   a_3}\Big|\\
\nnb &\le &\ds {12\L^2\over \l}\, \d+\Big| {  (b_3-b_1) (c_3-c_1)
\over
   a_3}-{(1-\d)^2 (b_3-b_1) (c_3-c_1)\over \d (1-\d)a_1+\d^2
a_2+(1-\d)^2a_3+\d(1-\d)a_4} \Big| \\
  &\le &\ds {12\L^2\over \l}\, \d+{8\L^3\over \l^2}\, \d \le
{20\L^3\over \l^2}\, \d.
\end{eqnarray}

Combining (\ref{E219}), (\ref{E233}) and (\ref{E235}), we get
(\ref{E216}) with $\ds C={56\L^3\sqrt \L\over \l^2 \sqrt \l }$.
\endpf

 The following result is concerned with the well-posedness and
regularity of state equation (\ref{E101}).

\begin{Lemma}\tt\label{T203}  Let {\rm(S1)--(S4)} hold and $z_0\in
L^\infty(\O)$. Then for any $u(\cd)\in\cU$, $(\ref{E101})$ admits a
unique weak solution $z(\cd)\in L^2(0,T;H^1_0(\O))\cap
L^{\infty}(\O_T)$. Furthermore, there exist a constant $K>0$,
independent of $u(\cd)\in\cU$, such that
\begin{equation}\label{E236A}
\|z(\cd)\|_{L^2(0,T;H^1_0(\O))}+\|z(\cd)\|_{L^\infty(\O_T)}\le
K.
\end{equation}
Moreover, there exists an $\a\in (0,1)$, such that  for any
$Q_0\subset\subset \O_T$, it holds that
\begin{equation}\label{E236}
 \|z(\cd) \|_{C^\a(Q_0)}\leq C(Q_0)
 \end{equation}
for some constant $C(Q_0)$.
\end{Lemma}
\Proof  The result is quite standard. We give a sketch of the proof.
Fix $u(\cd)\in \cU$. Let $m>0$, define
$$
f^m(t,x,z,u)=\left\{\begin{array}{ll} \ds f(t,x,z,u), & |z|\leq m, \\
\ds f(t,x,-m,u), &  z \leq -m, \\
\ds f(t,x,m,u), &  z \geq m.
\end{array}\right.
$$
For fixed $z(\cd)\in L^2(\O_T)$, let $z^m(\cd)$ be the solution of
$$
\left\{\begin{array}{ll} \ds\pa _t
z^m(t,x)-\div\big(A(t,x,u(t,x))\nabla
z^m(t,x)\big)=f^m(t,x,z(t,x),u(t,x)),\qq\eqin   \O_T,\\
\ds z^m(t,x)=0,  \qq\eqon [0,T]\times \pa\O,\\
\ds z^m(0,x)=z_0(x),\qq\eqin  \O,
\end{array}\right.
$$
Then there exist a constant $C_m>0$  such that
\begin{equation}\label{FE240}
 \|z^m(\cd)\|_{L^2(0,T;H^1_0(\O))}
+\|z^m(\cd)\|_{L^\infty(\O_T)}\leq C_m.
\end{equation}
Moreover, there exists an $\b=\b_m\in (0,1)$, such that  for any
$Q_0\subset\subset \O_T$, it holds that
\begin{equation}\label{FE241}
 \|z^m(\cd) \|_{C^{\b}(Q_0)}\leq C_m(Q_0)
\end{equation}
for some constant $C_m(Q_0)$.

Using (\ref{FE240})---(\ref{FE241}), we can see that the map
$z(\cd)\mapsto z^m(\cd)$ is continuous and compact from some ball of
$L^2(\O_T)$ to itself. Thus, Schauder fixed point theorem implies
that the map has a fixed point $Z^m(\cd)$. We have $Z^m(\cd)\in
L^2(0,T;H^1_0(\O))$ and
$$
\left\{\begin{array}{ll} \ns\ds\pa _t
Z^m(t,x)-\div\big(A(t,x,u(t,x))\nabla
Z^m(t,x)\big)=f^m(t,x,Z^m(t,x),u(t,x)),\qq\eqin   \O_T,\\
\ds Z^m(t,x)=0,  \qq\eqon  [0,T]\times \pa\O,\\
 \ds Z^m(0,x)=z_0(x),\qq\eqin  \O. \end{array}\right.
$$
Noting that (S4) holds, we can modify the proof of Theorem 7.1 of
Ch. 3  in  \cite{LSU} to  get that
$$
 \|Z^m(\cd)\|_{L^\infty(\O_T)}\leq C
$$
with $C$ being independent of $m$. Let $m>C$, we see that
(\ref{E101})  admits a unique weak solution $z(\cd)\in
L^2(0,T;H^1_0(\O))\cap L^{\infty}(\O_T)$ and (\ref{E236A}) holds.
Finally,  by (\ref{E106}),
\begin{equation}\label{E303}
|f(t,x,z (t,x),u (t,x))|\le M_K.
\end{equation}
Thus, (\ref {E236}) follows from Theorem 10.1 of Ch. 3  in
\cite{LSU}.
\endpf

\begin{Lemma}\tt\label{T204} Let $\d\in (0,1)$, $r>0$,
\begin{equation}\label{} \zeta_m(t,x)   =\left\{\begin{array}{ll}\ds
\d_{m1},& \eqif  (\{s\},\{x_1\})\in [\d,1)\times [\d,1),\\
\d_{m2},& \eqif   (\{s\},\{x_1\})\in [0,\d)\times [\d,1),\\
\d_{m3},& \eqif   (\{s\},\{x_1\})\in [\d,1)\times [0,\d),\\
\d_{m4},& \eqif   (\{s\},\{x_1\})\in [0,\d)\times
[0,\d),\end{array}\right. \q (t,x)\in \O_T.
\end{equation}
Then $\zeta_m({t\over \e^r},{x_1\over \e})$ ($m=1,2,3,4$) converges
weakly to $\mu_m$ in $L^2(\O_T)$ with
$$
\mu_1=(1-\d)^2, \q \mu_2=\mu_3=\d (1-\d), \q \mu_4=\d^2.
$$
\end{Lemma}
\Proof Such results are quite well-known and can be proved by
modifying the proof of Riemann's Lemma. One can verify easily that
for any rectangle $F  \subset \O_T$,
$$
\lim_{\e\to 0^+}\int_{\O_T} \chi_F(t,x) \zeta_1({t\over
\e^r},{x_1\over \e}) \, dt\, dx  =\mu_1\int_{\O_T} \chi_F(t,x) \,
dt\, dx.
$$
Since the set of all linear combinations of characteristic functions
$\chi_F(\cd)$ is dense in $L^2(\O_T)$, we get that $\zeta_1({t\over
\e^r},{x_1\over \e})$ converges weakly to $\mu_1$ in $L^2(\O_T)$.
The remains are similar.
\endpf

\begin{Lemma}\tt\label{T205} Let $n\ge 2$. Let $\xi,\eta\in\dbR^n$ be two
nonzero vectors. Then
\begin{equation}\label{1+xi*eta}
\sup_{|x|=1}\xi^\top
x x^\top \eta={|\xi|\,|\eta|+\xi^\top\eta\over 2},
\end{equation}
where ~$E^\top$ denotes the  transpose of a matrix ~$E$.
\end{Lemma}

The proof of above lemma is easy. See  \cite{Lou-Yong 2009}, for
example.

\vspace{6mm}

\def\theequation{3.\arabic{equation}}
\setcounter{equation}{0} \setcounter{section}{3}
\setcounter{Definition}{0} \setcounter{Remark}{0}\textbf{3. Proof of
the Main Theorem.} In this section, we present a proof of our main
theorem. The proof is divided into several steps. Let $\bar
u(\cd)\in\cU$ be an optimal control and $\bar z(\cd)$ be the
corresponding optimal state. Let $r>0$,
$u_2(\cd),u_3(\cd),u_4(\cd)\in\cU$ be fixed.

\textbf{I. Homogenizing spike variation of the control.}  Let $\d\in
(0,1)$ and $\e>0$.  For any $(t,x)=(t,x_1,x_2,\cds,x_n)\in \O_T$,
define
\begin{equation}\label{FE301} u^{\d,\e}(t,x)=\left\{\begin{array}{ll}\ds \bar u (t,x),& \eqif  (\{{t\over
\e^r}\},\{{x_1\over\e}\})\in [\d,1)\times [\d,1),\\
u_2(t,x),& \eqif   (\{{t\over
\e^r}\},\{{x_1\over\e}\})\in [0,\d)\times [\d,1),\\
u_3(t,x),& \eqif   (\{{t\over
\e^r}\},\{{x_1\over\e}\})\in [\d,1)\times [0,\d),\\
u_4(t,x),& \eqif   (\{{t\over \e^r}\},\{{x_1\over\e}\})\in
[0,\d)\times [0,\d).\end{array}\right. \end{equation}
Then $u^{\d,\e}(\cd)\in \cU$. Let $z^{\d,\e}(\cd)$ be the state
corresponding to $u^{\d,\e}(\cd)$, i.e.,
\begin{equation}\label{E301}\left\{\begin{array}{ll}
\ns\ds\pa _t z^{\d,\e}(t,x)-\div\big(A (t,x,u^{\d,\e}(t,x))\nabla
z^{\d,\e}(t,x)\big)=f(t,x,z^{\d,\e}(t,x),u^{\d,\e}(t,x)),\qq\eqin   \O_T,\\
\ns\ds z^{\d,\e}(t,x)=0,  \qq\hb{on
} [0,T]\times \pa\O,\\
\ns \ds z^{\d,\e}(0,x)=z_0(x),\qq\eqin  \O.
\end{array}\right.\end{equation}
By Lemma 2.3, there exists constants $K>0$ and $\a\in (0,1)$,
independent of $\d,\e$, such that
\begin{equation}\label{E302}\|z^{\d,\e}(\cd)\|_{L^2(0,T;H^1_0(\O))}+\|z^{\d,\e}(\cd)\|_{L^\infty(\O_T)}\le
K \end{equation}
and
\begin{equation}\label{E304B}
 \|z^{\d,\e}\|_{C^\a(Q_0)}\leq C(Q_0)  \end{equation}
for any  $ Q_0\subset\subset \O_T$ with some constant $C(Q_0)$.

By (\ref{E302}), for fixed $\d\in (0,1)$, we can extract a
subsequence (still denoted by itself) such that $z^{\d,\e}(\cd)$
converges to a function $z^\d(\cd)$ weakly in $ L^2(0,T;H^1_0(\O))$
as $\e\to 0^+$. By (\ref{E304B}) and Arzel\'a-Ascoli's theorem,
$z^{\d,\e}(\cd)$ converges uniformly to $z^\d(\cd)$ in $ C (Q_0) $
for any $Q_0\subset\subset \O_T$. Then, it follows easily from
$$
\|z^{\d,\e}(\cd)\|_{L^\infty(\O_T)}\le K
$$
that $z^{\d,\e}(\cd)$ converges  strongly to   $z^\d(\cd)$ in
$L^2(\O_T)$  and almost everywhere in $\O_T$.

By (\ref{E106}) and (\ref{E302}),
\begin{equation}\label{E304a} \Big|
f(t,x,z^{\d,\e}(t,x),u^{\d,\e}(t,x))-f(t,x,z^\d(t,x),u^{\d,\e}(t,x))\Big|
\le M_K|z^{\d,\e}(t,x)-z^\d(x)| .\end{equation}
On the other hand, by Lemma 2.4,  for any $h\in L^2(\O_T)$, when
$\e\to 0^+$,
\begin{eqnarray*}\label{E304C} \nnb & & \ds \int_{\O_T}   f(t,x,z^\d(t,x),u^{\d,\e}(t,x))h(t,x)\,dt\, dx \\
 \nnb &  \to &  \ds (1-\d)^2 \int_{\O_T} f(t,x,z^\d(t,x),\bar u(t,x))h(t,x)\,dt\, dx
 +\d (1-\d)  \int_{\O_T} f(t,x,z^\d(t,x), u_2(t,x))h(t,x)\,dt\, dx\\
& & \ds + \d (1-\d)  \int_{\O_T} f(t,x,z^\d(t,x),
u_3(t,x))h(t,x)\,dt\, dx+ \d^2 \int_{\O_T} f(t,x,z^\d(t,x),
u_4(t,x))h(t,x)\,dt\, dx.
\end{eqnarray*}
Combing the above with (\ref{E304a}), we get that along a
subsequence $\e\to 0^+$,
\begin{eqnarray*}\label{E304D}  \nnb & & \ds     f(t,x,z^{\d,\e}(t,x),u^{\d,\e}(t,x))  \\
 \nnb &  \to &  \ds (1-\d)^2   f(t,x,z^\d(t,x),\bar u(t,x))
 +\d (1-\d)    f(t,x,z^\d(t,x), u_2(t,x)) \\
  & & \ds + \d (1-\d)    f(t,x,z^\d(t,x), u_3(t,x)) + \d^2
  f(t,x,z^\d(t,x), u_4(t,x))  ,\q\mbox{weakly
in }\, L^2(\O_T).
\end{eqnarray*}
Let
\begin{eqnarray*}  h^{\d,\e}(t,x) &=  & \ds     f(t,x,z^{\d,\e}(t,x),u^{\d,\e}(t,x)) -(1-\d)^2   f(t,x,z^\d(t,x),\bar u(t,x)) \\
 & &  \ds
 -\d (1-\d)    f(t,x,z^\d(t,x), u_2(t,x)) - \d (1-\d)    f(t,x,z^\d(t,x), u_3(t,x))h(t,x) \\
&& \ds - \d^2
  f(t,x,z^\d(t,x), u_4(t,x))
\end{eqnarray*}
and $\tilde z^{\d,\e}(\cd)$ be the solution of
\begin{equation}\label{E305} \left\{\begin{array}{ll}
\ds\pa _t \tilde z^{\d,\e}(t,x)-\div\big(A
(t,x,u^{\d,\e}(t,x))\nabla \tilde
z^{\d,\e}(t,x)\big)=h^{\d,\e}(t,x),\qq\eqin   \O_T,\\
\ds \tilde z^{\d,\e}(t,x)=0,  \qq\eqon  [0,T]\times \pa\O,\\
\ds \tilde z^{\d,\e}(0,x)=0,\qq\eqin  \O.
\end{array}\right.
\end{equation}
Then
\begin{equation}\label{FE308}
 \int_\O\big|\tilde z^{\d,\e}(T,x)\big|^2\,
dx+\l\int_{\O_T}\big|\nabla \tilde z^{\d,\e}(t,x)\big|^2\, dt\,dx
\leq   \int_{\O_T} \tilde z^{\d,\e}(t,x)h^{\d,\e}(t,x)\,dt\,
dx.
\end{equation}
As  $z^{\d,\e}(\cd)$ converges  strongly in $L^2(\O_T)$, $\tilde
z^{\d,\e}(\cd)$ converges  strongly  in $L^2(\O_T)$ too.
Consequently,   it follows from \eqref{FE308} that
\begin{equation}\label{E306}
 \tilde z^{\d,\e}(\cd)\to 0,\qq\hb{strongly in }\,
L^2(0,T;H^1_0(\O)).
\end{equation}
By Lemma 2.1,
\begin{equation}\label{E307}  z^{\d,\e}(\cd)-\tilde  z^{\d,\e}(\cd) \to z^\d(\cd),\qq
\hb{weakly in}\, L^2(0,T;H^1_0(\O))
\end{equation}
with $z^\d(\cd)$ being the weak solution of
\begin{equation}\label{FE311}\left\{\begin{array}{ll}
\ds\pa_t z^\d (t,x)-\div\big(Q^\d (t,x)\nabla z^\d (t,x)\big)=&\ds
(1-\d)^2   f(t,x,z^\d(t,x),\bar u(t,x))\\
&\ds +\d (1-\d)    f(t,x,z^\d(t,x), u_2(t,x))  \\
&\ds  + \d (1-\d)
f(t,x,z^\d(t,x), u_3(t,x))\\
&\ds + \d^2
  f(t,x,z^\d(t,x), u_4(t,x)) ,\qq\eqin  \O_T,\\
\ds z^\d (t,x)=0,  &\qq\qq\eqon [0,T]\times \pa\O,\\
\ds z ^\d (0,x)=z_0(x),& \qq\qq\eqin  \O,
\end{array}\right.
\end{equation}
where $Q^\d(\cd)=\big(q^\d_{ij}(\cd)\big)\in L^\infty(\O_T;\cS^n_+)$
is given by
\begin{equation}\label{E312A}
q^\d_{ij}(t,x)=\int^1_0dy_1\int^1_0
\big(g^\d_{ij}(t,x,s,y_1)+
g^\d_{i1}(t,x,s,y_1)\pa_{y_1}\f^j(t,x,s,y_1)\big)ds \,,\q 1\le
i,j\le n
\end{equation}
and
\begin{equation}\label{}
\big(g^\d_{ij}(t,x,s,y_1)\big)=\left\{\begin{array}{ll}\ds
a_{ij} (t,x,\bar u(t,x)),& \eqif  (\{s\},\{y_1\})\in [\d,1)\times [\d,1),\\
\ds a_{ij}  (t,x, u_2(t,x)),& \eqif   (\{s\},\{y_1\})\in [0,\d)\times [\d,1),\\
\ds  a_{ij}  (t,x, u_3(t,x)),& \eqif   (\{s\},\{y_1\})\in [\d,1)\times [0,\d),\\
\ds  a_{ij}  (t,x, u_4(t,x)),& \eqif   (\{s\},\{y_1\})\in
[0,\d)\times [0,\d).\end{array}\right.
\end{equation}
For $r<2$,  $\f^k(t,x,\cd)\in L_\#^2(0,1;W^{1,2}_\#(0,1)/\dbR)$ is
the unique solution of
\begin{equation}\label{E206}
\pa_{y_1}
\big(g^\d_{1k}(t,x,s,y_1)+g^\d_{11}(t,x,s,y_1)\pa_{y_1} \f^k
(t,x,s,y_1)\big)=0.
\end{equation}
For $r=2$,  $\f^k(t,x,\cd)$ is the solution of
\begin{equation}\label{E207}
\pa_s \f^k(t,x,s,y_1)-\pa_{y_1}
\big(g^\d_{1k}(t,x,s,y_1)+g^\d_{11}(t,x,s,y_1)\pa_{y_1} \f^k
(t,x,s,y_1)\big)=0.
\end{equation}
While for $r>2$, $\f^k(t,x,s,y)=\f^k(t,x,y_1)$ is the solution of
\begin{equation}\label{E208BV}
 \pa_{y_1} \Big(\int^1_0g^\d_{1k}(t,x,s,y_1)\, ds
+\int^1_0 g^\d_{ij}(t,x,s,y_1)\, ds \, \pa_{y_1} \f^k
(t,x,y_1)\Big)=0.
\end{equation}
Combining (\ref{E307}) with (\ref{E306}), along a subsequence, we
obtain
\begin{equation}\label{E310}
z^{\d,\e}(\cd)\to z^\d(\cd),\qq \hb{weakly in
}\,L^2(0,T;H^1_0(\O)).
\end{equation}
Note that for fixed $\d\in (0,1)$, since any subsequence of
$z^{\d,\e}(\cd)$ has a further subsequence converging to the same
$z^\d(\cd)$ weakly in $L^2(0,T;H^1_0(\O))$, $z^{\d,\e}(\cd)$ itself
must converge to $z^\d(\cd)$ weakly in $L^2(0,T;H^1_0(\O))$.

In addition, by the optimality of $\bar u(\cd)$, we have
\begin{eqnarray}\label{E312}
\nnb \ds J(\bar u(\cd)) &\le & \lim_{\e\to
0^+}J\big(u^{\d,\e}(\cd)\big)\equiv J^\d =  \int_{\O_T}\Big((1-\d)^2
f^0(t,x,z^\d(t,x),\bar u(t,x)) \\
\nnb && \ds  +\d (1-\d)    f^0(t,x,z^\d(t,x), u_2(t,x)) \  + \d
(1-\d)
f^0(t,x,z^\d(t,x), u_3(t,x)) \\
 &&+ \d^2
  f^0(t,x,z^\d(t,x), u_4(t,x))\Big)\, dt\, dx.
\end{eqnarray}

\textbf{II. Linearized state equation.}  We now would like to let
$\d\to 0^+$. Denote
$$
Z^\d(x)={z^\d(x)-\bar z(x)\over\d},\qq (t,x)\in \O_T.
$$
Then it follows from (\ref{FE311}) that
\begin{equation}\label{E313}\left\{\begin{array}{l}
\ds\pa_t Z^\d (t,x) \ds -\div\big(A(t,x,\bar u(t,x))\nabla Z^\d
(t,x)\big)=\div\big({Q^\d(t,x)-A(t,x,\bar u(t,x))\over \d}\nabla
z^\d (t,x)\big)\\
\qq\qq\qq\ds +
(1-\d)^2   \int^1_0 f_z(t,x,\bar z(t,x)+s(z^\d(t,x)-\bar z(t,x)),\bar u(t,x))\, ds Z^\d (t,x)\\
\qq\qq\qq\ds +  (1-\d) \Big( f (t,x,z^\d(t,x),  u_2(t,x))-f (t,x,\bar z(t,x),  \bar u(t,x))\Big)\\
\qq\qq\qq\ds  + (1-\d) \Big( f (t,x,z^\d(t,x),  u_3(t,x))-f (t,x,\bar z(t,x),  \bar u(t,x))\Big)\\
\qq\qq\qq\ds + \d
   \Big( f (t,x,z^\d(t,x),  u_4(t,x))-f (t,x,\bar z(t,x),  \bar u(t,x))\Big) ,\qq\eqin  \O_T,\\
\ds Z^\d (t,x)=0,   \qq\qq\eqon [0,T]\times \pa\O,\\
 \ds Z ^\d (0,x)=0,  \qq\qq\eqin  \O.
\end{array}\right.\end{equation}
By (S3),
$$
|a_{ij}(t,x,v)|\le \L,\qq \forall\, (t,x)\in \O_T; v\in U;1\le
i,j\le n.
$$
Thus, it follows from (\ref{E312A}) and Lemmas 2.2 that as $\d\to
0^+$,
\begin{equation}\label{E315A}
{q^\d_{ij}(t,x)-a_{ij}(t,x,\bar u(t,x))\over\d}
\end{equation}
converges in $L^\infty(\O)$ to
\begin{eqnarray}\label{E315}
\nnb && \ds\th_{ij}(t,x) =
a_{ij}(t,x,u_2(t,x))+a_{ij}(t,x,u_3(t,x))-2a_{ij}(t,x,\bar u(t,x))\\
\nnb &&\ds -{\big[a_{1i}(t,x,u_3(t,x))-a_{1i}(t,x,\bar
u(t,x))\big]\big[a_{1j}(t,x,u_3(t,x))-a_{1j}(t,x,\bar
u(t,x))\big]\over a_{11}(t,x,u_3(t,x))}\,,\\
& &\hspace{7cm} \q 1\le i,j\le n.
\end{eqnarray}
On the other hand, $z^\d(\cd)$ is bounded uniformly in
$L^2(0,T;H^1_0(\O))$. Thus, we can prove step-by-step that as $\d\to
0^+$, $Z^\d(\cd)$ is bounded uniformly in $L^2(0,T;H^1_0(\O))$,
$z^\d(\cd)$ converges to $\bar z(\cd)$ strongly in
$L^2(0,T;H^1_0(\O))$, and $Z^\d(\cd)$ converges to $Z(\cd)$ weakly
in $L^2(0,T;H^1_0(\O))$ with $Z(\cd)$ being the weak solution of
\begin{equation}\label{E314}\left\{\begin{array}{l}
\ds\pa_t Z (t,x)  -\div\big(A(t,x,\bar u(t,x))\nabla Z
(t,x)\big)=\div\big(\Th(t,x)\nabla
\bar z  (t,x)\big)\\
\qq\qq\qq\ds +
  f_z(t,x,\bar z(t,x),\bar u(t,x)) \,  Z  (t,x)\\
\qq\qq\qq\ds +     f (t,x,\bar z(t,x),  u_2(t,x))+   f (t,x,\bar z(t,x),  u_3(t,x)) \\
\qq\qq\qq\ds  -2f (t,x,\bar z(t,x),  \bar u(t,x))   ,\qq\eqin  \O_T,\\
\ns\ds Z (t,x)=0,   \qq\qq\hb{on
} [0,T]\times \pa\O,\\
\ns \ds Z(0,x)=0,  \qq\qq\eqin  \O,
\end{array}\right.
\end{equation}
where
\begin{equation}\label{E314AA}\begin{array}{l}\ds
\Th(t,x)=\big(\th_{ij}(t,x)\big)=A(t,x,u_2(t,x))+A(t,x,u_3(t,x))-2A(t,x,\bar
u(t,x))\\
\, \\
\ds\qq -{\big[A(t,x,u_3(t,x))-A(t,x,\bar
u(t,x))\big]e_1e_1^\top\big[A(t,x,u_3(t,x))-A(t,x,\bar
u(t,x))\big]\over e_1^\top
A(t,x,u_3(t,x))e_1}\,.\end{array}
\end{equation}

From (\ref{E312}), we   have
\begin{eqnarray}\label{E315B}
\nnb  0 &\leq & \ds\lim_{\d\to 0^+}{J^\d-J(\bar u(\cd))\over \d} =\int_{\O_T}  \Big( f^0_z(t,x,\bar z(t,x),\bar u(t,x)) \,  Z  (t,x)\\
\nnb &&\ds +    f^0 (t,x,\bar z(t,x),  u_2(t,x)) + f ^0(t,x,\bar z(t,x),  u_3(t,x))\\
  && \ds -2f^0 (t,x,\bar z(t,x),  \bar u(t,x))\Big) \, dt\,
dx.
\end{eqnarray}

\textbf{III. Duality.} Let $\bar\psi(\cd)$ be the solution of the
adjoint equation (\ref{adjoint}). Then (\ref{E315B}) becomes
\begin{eqnarray}\label{E317}
\nnb 0 &\le & \int_{\O_T}\Big[f^0 (t,x,\bar z(t,x),  u_2(t,x)) + f ^0(t,x,\bar z(t,x),  u_3(t,x))-2f^0 (t,x,\bar z(t,x),  \bar u(t,x)) \\
\nnb && + \Big(\pa_t \psi (t,x)+\div\big(A(t,x,\bar u(t,x))\nabla
\bar\psi
(t,x)\big)+f_z(t,x,\bar z(t,x),\bar u(t,x))\,\bar\psi(t,x)\Big) Z(t,x)\Big]\,dt\,dx\\
\nnb &=& \ds
\int_{\O_T}\Big[f^0 (t,x,\bar z(t,x),  u_2(t,x)) + f ^0(t,x,\bar z(t,x),  u_3(t,x))-2f^0 (t,x,\bar z(t,x),  \bar u(t,x)) \\
\nnb && \ds +\Big(-\pa_t Z (t,x)+\div\big(A(t,x,\bar u(t,x))\nabla
Z(t,x)\big)+f_z(t,x,\bar z(t,x),\bar u(t,x))\,Z(t,x)\Big)\bar\psi(t,x)\Big]\,dt\,dx
\\
\nnb &=& \ds \int_{\O_T}\Big[ f^0 (t,x,\bar z(t,x),  u_2(t,x)) + f ^0(t,x,\bar z(t,x),  u_3(t,x))-2f^0 (t,x,\bar z(t,x),  \bar u(t,x)) \\
\nnb && \ds - \Big(f (t,x,\bar z(t,x),  u_2(t,x)) + f (t,x,\bar
z(t,x), u_3(t,x))-2f (t,x,\bar z(t,x),  \bar u(t,x))\Big)\,
\bar\psi(t,x) \\
\nnb && + \lan\Th(t,x)\nabla\bar z(t,x),\nabla\bar\psi(t,x)\ran\,\Big] \, dt\, dx\\
\nnb &=& \ds \int_{\O_T}\Big[2\( f (t,x,\bar z(t,x),  \bar
u(t,x))\bar\psi(x)- f^0 (t,x,\bar z(t,x),  \bar u(t,x))- \lan
A(t,x,\bar u(t,x))\nabla\bar z(t,x),\nabla\bar\psi(t,x)\ran\)
 \\
\nnb && \ds -\( f (t,x,\bar z(t,x),   u_2(t,x))\bar\psi(x)- f^0
(t,x,\bar z(t,x),  u_2(t,x))-\lan
A(t,x,u_2(t,x))\nabla\bar z(t,x),\nabla\bar\psi(t,x)\ran\)  \\
\nnb && \ds  -\( f (t,x,\bar z(t,x),   u_3(t,x))\bar\psi(x)- f^0
(t,x,\bar z(t,x),  u_3(t,x))-\lan
A(t,x,u_3(t,x))\nabla\bar z(t,x),\nabla\bar\psi(t,x)\ran\)\\
\nnb && \ds -{\nabla\bar z(t,x)^\top \big(A(t,x,u_3(t,x))-A(t,x,\bar
u(t,x))\big)e_1e_1^\top\big(A(t,x,u_3(t,x))-A(t,x,\bar
u(t,x))\big)\nabla\bar\psi(t,x)\over e_1^\top A(t,x,u_3(t,x))e_1}\Big]\, dt\, dx\\
\nnb &=& \ds\int_{\O_T}\Big[2H\big(t,x,\bar
z(t,x),\bar\psi(t,x),\nabla\bar
z(t,x),\nabla\bar\psi(t,x),\bar u(t,x)\big) \\
\nnb &&  -H\big(t,x,\bar z(t,x),\bar\psi(t,x),\nabla\bar
z(t,x),\nabla\bar\psi(t,x),
u_2(t,x)\big) \\
\nnb && \ds -H\big(t,x,\bar z(t,x),\bar\psi(t,x),\nabla\bar
z(t,x),\nabla\bar\psi(t,x), u_3(t,x)\big) \\
 && \ds - \F\big(A(t,x,u_3(t,x))-A(t,x,\bar
u(t,x)),A(t,x,u_3(t,x)),\nabla\bar z(t,x),\nabla
\bar\psi(t,x),e_1\big)\Big] \,dt\, dx,
\end{eqnarray}
where $H$ is defined by (\ref{E110}),
\begin{equation}\label{E319}
\F\big(A,B,\xi,\eta,\n\big)={\eta^\top A\n\n^\top A \xi\over \n^\top
B\n}\,,\qq (A,B,\xi,\eta,\n)\in\cS^n\times\cS^n_+\times\dbR^n\times
\dbR^n\times\dbR^n
\end{equation}
and $\cS^n$ is the set of all $n\times n$ real symmetric matrices.

\textbf{IV. Maximum condition.} By a standard argument
(\cite{Li-Yong 1995}), it follows from (\ref{E317}) that
\begin{eqnarray}\label{E320}
\nnb &  & \ds 2H\big(t,x,\bar z(t,x),\bar\psi(t,x),\nabla\bar
z(t,x),\nabla\bar\psi(t,x),\bar u(t,x)\big) \\
\nnb &  \geq &\ds  H\big(t,x,\bar z(t,x),\bar\psi(t,x),\nabla\bar
z(t,x),\nabla\bar\psi(t,x),
w \big) \\
\nnb &  & \ds +H\big(t,x,\bar z(t,x),\bar\psi(t,x),\nabla\bar
z(t,x),\nabla\bar\psi(t,x), v\big) \\
\nnb &  & \ds + \F\big(A(t,x,v)-A(t,x,\bar
u(t,x)),A(t,x,v),\nabla\bar
z(t,x),\nabla \bar\psi(t,x),e_1\big),  \\
 &  & \qq\qq\qq\forall v,w\in U,\q\ae (t,x)\in\O_T.
\end{eqnarray}
Further, it is easy to see that (\ref{E320}) is equivalent to the
following two inequalities:
\begin{eqnarray} \label{E320B}
\nnb &  & \ds  H\big(t,x,\bar z(t,x),\bar\psi(t,x),\nabla\bar
z(t,x),\nabla\bar\psi(t,x),\bar u(t,x)\big) \\
 & \geq &\ds   H\big(t,x,\bar z(t,x),\bar\psi(t,x),\nabla\bar
z(t,x),\nabla\bar\psi(t,x), v \big),  \qq\forall v\in U,\q\ae
(t,x)\in\O_T
\end{eqnarray}
and
\begin{eqnarray} \label{E320C}
\nnb &  &  H\big(t,x,\bar z(t,x),\bar\psi(t,x),\nabla\bar
z(t,x),\nabla\bar\psi(t,x),\bar u(t,x)\big)\\
\nnb &  \geq  & \ds  H\big(t,x,\bar z(t,x),\bar\psi(t,x),\nabla\bar
z(t,x),\nabla\bar\psi(t,x), v\big) \\
\nnb &  & \ds + \F\big(A(t,x,v)-A(t,x,\bar
u(t,x)),A(t,x,v),\nabla\bar
z(t,x),\nabla \bar\psi(t,x),e_1\big),  \\
  &  & \qq\qq\qq\forall v \in U,\q\ae (t,x)\in\O_T.
\end{eqnarray}
Moreover, we   can generalize (\ref{E320C}) to the following:
\begin{eqnarray}\label{E320CC}
\nnb & & \ds  H\big(t,x,\bar z(t,x),\bar\psi(t,x),\nabla\bar
z(t,x),\nabla\bar\psi(t,x),\bar u(t,x)\big)\\
\nnb &  \geq  & \ds  H\big(t,x,\bar z(t,x),\bar\psi(t,x),\nabla\bar
z(t,x),\nabla\bar\psi(t,x), v\big) \\
\nnb &  & \ds + \F\big(A(t,x,v)-A(t,x,\bar
u(t,x)),A(t,x,v),\nabla\bar
z(t,x),\nabla \bar\psi(t,x),e \big),  \\
\nnb &  & \qq\qq\qq\forall v \in U,e\in S^{n-1}, \q\ae (t,x)\in\O_T,
\end{eqnarray}
where $S^{n-1}=\{x\in\dbR^n\bigm||x|=1\}$. Now, for given $(t,x,v)$,
we denote
\begin{equation}\label{}\left\{\begin{array}{ll}
\ds\m={A(t,x,v)^{1\over2}  e\over |A(t,x,v)^{1\over2}  e|} ,\\
\ds\xi=A(t,x,v)^{-{1\over 2}}\big[A(t,x,v)-A(t,x,\bar
u(t,x))\big]\nabla\bar z(t,x),\\
\ds\eta=A(t,x,v)^{-{1\over2}}\big[A(t,x,v)-A(t,x,\bar
u(t,x))\big]\nabla\bar\psi(t,x).\end{array}\right.
\end{equation}
When $e$ runs over $S^{n-1}$, $ \m$ will run over $S^{n-1}$. Then (\ref{E320B}) and (\ref{E320C}) become
\begin{eqnarray}\label{E325}
\nnb & &  \ds H\big(t,x,\bar z(t,x),\bar\psi(t,x),\nabla\bar
z(t,x),\nabla\bar\psi(t,x),\bar u(t,x)\big)\\
\nnb &   & \ds - H\big(t,x,\bar z(t,x),\bar\psi(t,x),\nabla\bar
z(t,x),\nabla\bar\psi(t,x), v\big) \\
  &   \geq & \ds \sup_{|\m|=1} \max \Big( \xi^\top \m\m^\top
\eta , 0\Big)\,.
\end{eqnarray}
By
Lemma \ref{T205}, one can get that
$$
\sup_{|\mu|=1} \xi^\top \m\m^\top\eta= \left\{\begin{array}{ll}\ds
{|\xi|\,|\eta|+\xi^\top\eta\over 2}\, & {\rm if}\, n\geq 2,\\
\xi\, \eta\,  & {\rm if}\, n=1.\end{array}\right.
$$
Thus
\begin{equation}\label{}
 \sup_{|\m|=1} \max \Big(  \xi^\top \m\m^\top\eta,
0\Big)={|\xi|\,|\eta|+\xi^\top\eta\over 2}\,.
\end{equation}
Combining the above, we obtain (\ref{E109}). This completes the
proof of Theorem \ref{T101}.
\endpf

We can see that the limit equation \eqref{FE311} of homogenizing
spike variation equation \eqref{E301} may be different for different
$r$. More precisely, \eqref{FE311} has essentially three different
cases corresponding to $r<2$, $r=2$ and $r>2$, respectively.
However, the variational equation \eqref{E313} is independent of
$r>0$. Thus, the final result (Theorem 1.1) can be got by choosing
$r\in (0,2)$. Such a choice will lead to a simple proof of Theorem
1.1. But, if we did that, we would not know whether we can get other
conditions from cases of $r\geq 2$. This is not satisfied.

Concerning the method to construct spike variation, we have
mentioned that special forms of spike variation are needed to get
good expressions of the limit equations as \eqref{FE311}. When we
introduce  \eqref{FE301}, it is natural to expect that $u_4(\cd)$
has no essential effect on the final result. The effect of
$u_2(\cd)$ is in time scale. One can see that essentially,
$u_2(\cd)$ works as a spike variation as $A(\cd)$ being independent
of $u$. Difficulties caused by $A(\cd)$ containing $u$ appear when
$u_3(\cd)$ is introduced.

\begin{Remark}\tt\label{R301}
If we follow the idea of sequential laminates (see Tartar
\cite{Tar}), we can generalize \refeq{FE311} and \refeq{E312} by
constructing more general homogenized equations with their leading
terms satisfying
\begin{eqnarray*}
  Q^\d(t,x) &=& \bA(t,x)+\d \big(A_2(t,x)+A_3(t,x)-2 \bA(t,x)\big)\\
  && -\d \big(A_3(t,x)-\bA(t,x)\big) \big(A_3(t,x)\big)^{-{1\over 2}}P\,
\big(A_3(t,x)\big)^{-{1\over 2}}\big(A_3(t,x)-\bA(t,x)\big)+o(\d),
\end{eqnarray*}
where
$$
\bA (t,x)=A(t,x,\bu(t,x)), \q A_2(t,x)=A(t,x,u_2(t,x)), \q A_3(t,x)=A(t,x,u_3(t,x))
$$
and
\begin{eqnarray*}
P\in \mcS_1 & \defeq &\{W\in \cS^n| \,  W\ge 0, \tr W=1\}\\
&=& \{\sum^n_{k=1}\a_k\xi_k\xi_k^\top | \a_1,\a_2,\ldots,\a_n\in
[0,1],\xi_1,\xi_2,\ldots,\xi_n\in S^{n-1}, \sum^n_{k=1}\a_k=1 \}.
\end{eqnarray*}
Consequently, we can generalize \refeq{E314}---\refeq{E315B} with \refeq{E314AA} being replaced by
\begin{eqnarray}\label{E332}
\nnb  \GT(t,x)&=& A_2(t,x)+A_3(t,x)-2 \bA(t,x)\\
  && -\big(A_3(t,x)-\bA(t,x)\big) \big(A_3(t,x)\big)^{-{1\over 2}}P\,
\big(A_3(t,x)\big)^{-{1\over 2}}\big(A_3(t,x)-\bA(t,x)\big).
\end{eqnarray}
It follows from
$$
\sup_{P\in \mcS_1} \xi^\top P \eta =\sup_{|x|=1}\xi^\top x x^\top \eta
$$
that \refeq{E314}, \refeq{E315B} and \refeq{E332} still lead to Theorem \ref{T101}.
\end{Remark}

\vspace{6mm}

\def\theequation{4.\arabic{equation}}
\setcounter{equation}{0} \setcounter{section}{4}
\setcounter{Definition}{0} \setcounter{Remark}{0}\textbf{4. Problem
with State Constraints.} In this section, we will consider the cases
of state constraint. We will only state the results since the proofs
are completely similar to those of  elliptic cases.

(S6)  Let $\cZ$ be a Banach space with strictly convex dual $\cZ^*$,
$F: L^2(0,T;H^1_0(\O))\to \cZ$ be  continuous Fr\'echet
differentiable, and $E\subseteq \cZ$ be closed and convex.

As in Chapter 5 of \cite{Li-Yong 1995}, many state constraints can
be stated in the following type:
\begin{equation}\label{E501}
F(z(\cd))\in E.
\end{equation}
Let $\cP_{ad}$ be the set of all pairs $(z(\cd),u(\cd))$ satisfying
\eqref{E101} and \eqref{E501}. Any $(z(\cd),u(\cd))\in \cP_{ad}$ is
called an admissible pair. The set $\cU_{ad}\equiv
\{u(\cd)\in\cU\,|\, (z(\cdot;u(\cdot)),u(\cdot))\in \cP_{ad}\}$ is
called the set of admissible controls. Our optimal control problem
with state constraint is

\textbf{Problem (SC)}. Find a control $\bar u(\cd)\in\cU_{ad}$ such
that
\begin{equation}\label{E502}
J(\bar u(\cd))=\inf_{ u(\cd) \in\cU}J(u(\cd)).
\end{equation}
To state necessary conditions for optimal controls of Problem (SC),
we need to recall the notion of finite co-dimensionality (see Ch. 4
of \cite{Li-Yong 1995}, for example).
\begin{Definition}\tt\label{T401}
Let $X$ be a Banach space and $X_0$ be a subspace of $X$. We say
that $X_0$ is finite co-dimensional in $X$ if there exist
$x_1,x_2,\cds,x_n\in X$, such that
\begin{equation}\label{E503}
\span\{X_0,x_1,\cds,x_n\}\equiv\hb{the space spanned by
$\{X_0,x_1,\cds,x_n\}$}=X.
\end{equation}
A subset $S$ of $X$ is said to be finite co-dimensional in $X$ if
for some $x_0\in S$, $\span (S-\{x_0\})\equiv$ the closed subspace
spanned by $\{x-x_0|x\in S\}$ is a finite co-dimensional subspace of
$X$ and $\coh S\equiv $ the closed convex hull of $S-\{x_0\}$ has a
nonempty interior in this subspace.
\end{Definition}

Let $(\bar z(\cd),\bar u(\cd))$ be an optimal pair of Problem (SC).
Let $Z=Z(\cd;u(\cd))\in L^2(0,T; H^1_0(\O))$ be the unique weak
solution of the variational equation \eqref{E314} and define the
reachable set of variational system \eqref{E314} as
\begin{equation}\label{varR}
\cR=\big\{Z(\cd; u(\cd))\bigm|u_2(\cd),
u_3(\cd)\in\cU\big\}.
\end{equation}

Now, let us state the necessary conditions of an optimal control to
Problem (SC) as follows:

\begin{Theorem}\tt\label{T402} Let {\rm(S1)--(S6)} hold. Let $(\bar z(\cd),\bar
u(\cd))\in \cP_{ad}$ be an optimal pair of Problem {\rm(SC)}. Let
$$
F'(\bar z(\cd))\cR-E\equiv\{\xi-\eta\bigm|\xi\in F'(\bar z(\cd))\cR,\eta\in E\}
$$
be finite co-dimensional in $\cZ$. Then there exists a triple
$(\bar\psi_0,\bar\psi(\cd),\bar\f(\cd))\in\dbR\times L^2(0,T;
H^1_0(\O))\times\cZ^*$ satisfying
\begin{equation}\label{nontr}
\left\{\begin{array}{ll}
\ds\bar\psi_0\le 0,\\
\ds(\bar\psi_0,\bar\f(\cd))\ne0,\\
\ds(\bar\psi_0,\bar\psi(\cd))\ne 0,\qq\hb{ if }F'(\bar z(\cd))^*
\mbox{ is injective},\end{array}\right.
\end{equation}
\begin{equation}\label{trans}\lan\bar\f(\cd),\eta-F(\bar
z(\cd))\ran{}_{\cZ^*\1n,\,\cZ}\le0, \qq\qq\forall\eta\in
E,
\end{equation}
\begin{equation}\label{adjointS}\left\{\begin{array}{ll}
 \ds \pa_t \bar \psi (t,x)+\div\big(A(t,x,\bar
u(t,x))\nabla\bar\psi(t,x)\big)=-\psi_0 f^0_z(t,x,\bar z(t,x),\bar
u(t,x))\\
 \ds \qq\qq -f_z(t,x,\bar z(t,x),\bar
u(t,x))\,\bar\psi(t,x)+F^\prime(\bar z(\cd))^*\bar\f
,\qq\eqin  \O_T,\\
\ds\bar\psi(t,x) =0,\qq\eqon [0,T]\times \pa\O,\\
\ds \bar\psi(T,x)=0, \qq\eqin \O,\end{array}\right.
\end{equation}
\begin{eqnarray}\label{E408}
\nnb & & \ds H\big(t,x,\bar z(t,x),\bar\psi(t,x),\nabla\bar
z(t,x),\nabla\bar\psi(t,x),\bar u(t,x)\big)
-H\big(t,x,\bar z(t,x),\bar\psi(t,x),\nabla\bar z(t,x),\nabla\bar\psi(t,x),v\big)\\
\nnb &  \ge &\ds {1\over2}\big|A(t,x,v)^{-{1\over 2}}(A(t,x,\bar
u(t,x))-A(t,x,v))\nabla\bar z(t,x)\big|\,\big|A(t,x,v)^{-{1\over
2}}(A(t,x,\bar u(t,x))-A(t,x,v))\nabla\bar\psi(t,x)\big|\\
\nnb &  & \ds +{1\over2}\lan A(t,x,v)^{-{1\over2}}(A(t,x,\bar
u(t,x))-A(t,x,v))\nabla\bar z(t,x),A(t,x,v)^{-{1\over2}}(A(t,x,\bar
u(t,x))-A(t,x,v))\nabla\bar\psi(t,x)\ran,\\
  &  &\ds \qq\qq\qq\qq\forall v\in U,\q\ae
(t,x)\in\O_T,
\end{eqnarray}
where
\begin{eqnarray}\label{E409}
\nnb &  &\ds H(t,x,z,\psi,\xi,\eta,v)=\lan\psi,f(t,x,z,v)\ran+\psi_0
f^0(t,x,z,v)
-\lan A(t,x,v)\xi,\eta\ran,\\
  & &\ds\qq\qq\qq \qq (t,x,z,\psi,\xi,\eta, v)\in
[0,T]\times \O\times\dbR\times \dbR\times\dbR^n\times\dbR^n\times
U.
\end{eqnarray}
\end{Theorem}

\textbf{Acknowledgement.} The author  thanks the anonymous referees
for their helpful suggestions.

 \footnotesize
\vspace{5mm}\footnotesize \ \\

\end{document}